\hsize 159.2mm
\vsize 246.2mm
\font\Bbb=msbm10

\font\bigrm=cmr17

\magnification=\magstep1
\def\C{\hbox{\Bbb C}}

\centerline{\bigrm A Decomposition of Complex Monge-Amp\`ere Measures }
\vskip .3in
\centerline{\sl Yang Xing     }
\vskip .6in
\noindent {{ Abstract.}  We prove one decomposition theorem of complex Monge-Amp\`ere measures of plurisubharmonic functions in connection with their pluripolar sets.
}\bigskip

\bigskip\bigskip 
\noindent{\bf 1. Introduction}
\bigskip

The purpose of this paper is to give a decomposition of complex Monge-Amp\`ere measures associated to pluripolar sets of plurisubharmonic functions in the class ${\cal F}(\Omega)$ defined in [C1]. We denote by $PSH(\Omega)$ the class of plurisubharmonic functions in a hyperconvex domain $\Omega$ and by $PSH^-(\Omega)$ the subclass of negative functions. Recall that a set $\Omega\subset \C^n$ is said to be a hyperconvex domain if it is open, bounded, connected and there exists $\rho\in PSH^-(\Omega)$ such that $\{z\in\Omega;\rho(z)<-c\}\subset\subset\Omega$ for any $c>0$.
The class ${\cal F}(\Omega)$ consists of all plurisubharmonic functions $u$ in $\Omega$ such that there exists a sequence $u_j\in {\cal E}_0(\Omega)$, $u_j\searrow u,\,j\to\infty$ and $\sup_j\int_\Omega(dd^cu_j)^n<\infty$, where ${\cal E}_0(\Omega)$ is the class of bounded plurisubharmonic functions $v$ with $\lim_{z\to\zeta}v(z)=0$ for all $\zeta\in \partial\Omega$ and $\int_\Omega(dd^cv)^n<\infty.$ We also need the subclass ${\cal F}^a(\Omega)$ of functions from ${\cal F}(\Omega)$ whose Monge-Amp\`ere measures put no mass on pluripolar subsets of $\Omega$. 
It is known that Monge-Amp\`ere measures $(dd^cu)^n$ for $u\in {\cal F}(\Omega)$  are  well-defined finite measures in $\Omega$, see [C1] for details.
Our main result is the following: Restriction of the complex Monge-Amp\`ere measure of a function $u\in {\cal F}(\Omega)$ on its pluripolar set is still a Monge-Amp\`ere measure of some function in ${\cal F}(\Omega)$. As an application we obtain that every Monge-Amp\`ere measure of  functions in ${\cal F}(\Omega)$ can be written as a sum of two Monge-Amp\`ere measures where one  has zero mass on any pluripolar set and another one is carried by the pluripolar set of the corresponding function.
\bigskip
It is a great pleasure for me to thank Urban Cegrell for many fruitful comments.
\bigskip\bigskip 
\noindent{\bf 2. Theorems and Proofs}
\bigskip
We need an inequality.
\bigskip
\noindent {\bf Lemma.[X2].} \it Let $u,\,v\in
PSH(\Omega)\,\cap\,L^\infty(\Omega)$ be such that $\liminf\limits_{z\to\partial\Omega} \bigl(
u(z)-v(z)\bigr)\geq 0.$ Then for any $-1\leq w\in PSH^-(\Omega)$ we have that $ (n!)^{-2} \int_{u<v} (v-u)^n\,(dd^cw)^n+\int_{u<v} (-w)\,(dd^cv)^n\leq \int_{u<v}
(-w)\,(dd^cu)^n. $
\rm \bigskip
Recall [X2] that a sequence $\{u_j\}$ of functions in $PSH(\Omega)$ is said to be convergent to a function $u$ in $C_n$ on a subset $E$ of $\Omega$ if  for any $\delta>0$ we have that $C_n\{z\in E;|u_j(z)-u(z)|>\delta\}\to 0$ as $ j\to\infty$, where $C_n$ denotes the inner capacity introduced by Bedford and Taylor in [BT].
\bigskip
\noindent {\bf Theorem 1.}  \it  Let $v\in {\cal F}(\Omega)$. Then there exists $u\in {\cal F}(\Omega)$ with $u\geq v$ in $\Omega$ such that
$$(dd^cu)^n=\chi_{\{v=-\infty\}}\,(dd^cv)^n \qquad {\rm in}\quad \Omega,$$ 
where $\chi_{\{v=-\infty\}}$ is the characteristic
function of $\{v=-\infty\}$. Furthermore, let $g$ be the unique function in ${\cal F}^a(\Omega)$ with $(dd^cg)^n =\chi_{\{v>-\infty\}}\,(dd^cv)^n$, then $v\geq u+g$ in $\Omega$.
\rm 
\bigskip
\noindent{\it Proof.} By Theorem 2.1 in [C1] we can take a sequence $v_j\in {\cal E}_0(\Omega)$ such that $v_j\searrow v$, $j\to\infty$. 
By [C2][K] there exist $u^k_j\in
{\cal E}_0(\Omega)$ such that 			
$(dd^cu^k_j)^n=-\max(v/k,\,-1)\,(dd^cv_j)^n$. 
From the comparison theorem [BT] it follows that $ u^{k+1}_j\geq u^k_j\geq v_j\geq v$.
By passing to a subsequence if necessary, we assume that $u^k_j\to u^k\in {\cal F}(\Omega)$ weakly, $j\to\infty$, and $u^k\nearrow u\in {\cal F}(\Omega)$, $k\to\infty$. Then by Theorem 2 below we have that $(dd^cu^k)^n=-\max(v/k,\,-1)\,(dd^cv)^n$, which implies  $(dd^cu)^n=\chi_{\{v=-\infty\}}\,(dd^cv)^n.$
If furthermore $\chi_{\{v>-\infty\}}\,(dd^cv)^n=(dd^cg)^n $ for $g\in {\cal F}^a(\Omega)$, then we take $g_j^k\in {\cal E}_0(\Omega)$ such that $(dd^cg_j^k)^n=\max((v+k)/k,\,0)\,(dd^cv_j)^n=\max((v+k)/k,\,0)\,
(dd^c\max (v_j,-k-1))^n.$ 
By the comparison theorem [BT] we have $0>g_j^k\geq \max (v_j,-k-1)\geq v$.
By Theorem 2 below again, we  assume that $g_j^k$ converges to a 
bounded psh function $g^k$ in $C_n$ on each $E\subset\subset \Omega$. 
Letting $j\to \infty$ we get that
$(dd^cg^k)^n=\max((v+k)/k,\,0)\,(dd^cv)^n=\max((v+k)/k,\,0)\,(dd^cg)^n\leq (dd^cg)^n$, which implies  
$0>g^k\geq g$. Hence $g^k$ decrease to some $g_1\in {\cal F}^a(\Omega)$. By Theorem 5.15 in [C1] we have $g_1=g$.
Since $\bigl(dd^c(g_j^k+u_j^k)^n\bigr)^n\geq (dd^cg_j^k)^n+(dd^cu_j^k)^n=
(dd^cv_j)^n$ we get that $v_j\geq g_j^k+u_j^k$ and hence $ v \geq g+u$. The proof of Theorem 1 is
complete.
\bigskip
\noindent {\bf Theorem 2.}  \it Suppose that $v\in {\cal F}(\Omega)$, $v_j\in {\cal E}_0(\Omega)$ and  $-1\leq \psi\in PSH^-(\Omega)$ are such that $v_j\searrow v$ as $j\to\infty$ and $v$ is bounded on $\{z\in\Omega;\psi(z)\not=-1\}$. If $u_j\in {\cal E}_0(\Omega)$ are such that $(dd^cu_j)^n=-\psi\,(dd^cv_j)^n$ and $u_j\to u\in PSH(\Omega)$ weakly in $\Omega$, then $(dd^cu)^n=-\psi\,(dd^cv)^n$, $u\geq v$ and hence $u\in{\cal F}(\Omega)$.
\rm 
\bigskip
\noindent{\it Proof.} Clearly, $0\geq u_j\geq v_j\geq v$. Hence $u\geq v$ and $u\in{\cal F}(\Omega)$.
To prove $(dd^cu)^n=-\psi\,(dd^cv)^n$, by Theorem 7 in [X1] or [C1] we have that $-\psi\,(dd^cv_j)^n\to -\psi\,(dd^cv)^n$ weakly, $j\to\infty$, and hence it is
enough to show that $u_j\to u $ in $C_n$ on each $E\subset\subset\Omega$ as $j\to\infty$.
Take $t<\inf_{\{\psi\not=-1\}}v$. Since $(dd^cv_j)^n=\chi_{\{v_j>t\}}\,(dd^cv_j)^n+\chi_{\{v_j\leq t\}}\,(dd^cv_j)^n\leq \bigl(dd^c\max(v_j,t)\bigr)^n+(dd^cu_j)^n\leq \Bigl(dd^c\bigl(\max(v_j,t)+u_j\bigr)\Bigr)^n$,
we have $v_j\geq u_j+\max(v_j,t)$ and thus $v\geq u+t$.
Given $E\subset\subset\Omega$ and $0<\varepsilon<-t$ take $0<\delta<1$ such that
$C_n\bigl\{z\in E;(1-\delta)\,v\leq -\varepsilon\bigr\}<\varepsilon$.
By quasicontinuity of psh functions and Hartog's Lemma, we only need to show that $C_n\bigl\{z\in E;\ 
u(z)>u_j(z)+3\varepsilon\bigr\}\longrightarrow 0$, $j\to\infty$. Let $l_j:=\min_\Omega \bigl(\delta\,u_j+\varepsilon\bigr)$. 
Since 
$C_n\bigl\{z\in E;\ u_j(z)\leq\delta \,u_j(z)-\varepsilon\bigr\}\leq  C_n\bigl\{z\in E;(1-\delta)\,v\leq -\varepsilon\bigr\}<\varepsilon$, 
we have $$C_n\bigl\{z\in E;\  u(z)>u_j(z)+3 \varepsilon\bigr\}\leq
C_n\bigl\{z\in\Omega;\  u(z)>\delta\,u_j(z)+2 \varepsilon\bigr\}+\varepsilon$$
$$\leq \sup \biggl\{ \
{1\over \varepsilon^n}\int\limits_{
u>\delta\,u_j+\varepsilon}\bigl(u-\delta\,u_j-\varepsilon\bigr)^n\,(dd^cw)^n;\ w\in PSH(\Omega),\
0<w<1\ \biggr\}+\varepsilon$$ $$= \sup \biggl\{ \ {1\over
\varepsilon^n}\ \int\limits_{ \max(u,l_j)\
>\delta\,u_j+\varepsilon}\bigl(\max(u,l_j)-\delta\,u_j-\varepsilon\bigr)^n(dd^cw)^n;$$  $$ w\in
PSH(\Omega), 0<w<1\ \biggr\}+\varepsilon$$ 
which by Lemma is less than
$${(n!)^2\delta^n\over
\varepsilon^n}\ \int\limits_{ \max( u,l_j)\
>\delta\,u_j+\varepsilon}(dd^c u_j)^n+\varepsilon\leq {(n!)^2\delta^n\over
\varepsilon^n}\int\limits_{ u
>\delta\,u_j+\varepsilon}(dd^c v_j)^n+\varepsilon$$
 $$\leq {(n!)^2\delta^n\over
\varepsilon^n}\int\limits_{ u
>\delta\,u_j+\varepsilon}\phi\,(dd^c v_j)^n+2\,\varepsilon$$
for some $\phi\in C_0^\infty(\Omega)$ with $0\leq \phi\leq 1$, where the last inequality follows from the fact that $\lim_{j\to\infty}\int_\Omega(dd^cv_j)^n=\int_\Omega(dd^cv)^n<\infty$.
Since $v-t\geq u\geq v$ and $ u_j\geq v$, then for $a=(\varepsilon+t)/(1-\delta)<0$ the last integral equals
$$ \int\limits_{{\rm max} (u,a)>\delta\,{\rm max} (u_j,a) +\varepsilon}\phi\,(dd^cv_j)^n$$
$$\leq {1\over
\varepsilon}\,\int\limits_{
{\rm max} (u,a)>\delta\,{\rm max} (u_j,a) +\varepsilon}\phi\,\bigl({\rm max} (u,a)-{\rm max} (u_j,a) \bigr)   
(dd^cv_j)^n.$$
By $v_j\geq u_j+t$ and  $v_j\geq v\geq
u+t$ we have that ${\rm max} (u,a)-{\rm max} (u_j,a)=0$ if  $v_j\leq a+t$. 
By quasicontinuity of the $u$ there exists an open subset $O_\varepsilon\subset\Omega$ such that $C_n(O_\varepsilon)<\varepsilon^{n+2}$ and $u\in C(\Omega\setminus O_\varepsilon)$.
It then follows from Hartog's Lemma  that $\varepsilon^{n+2}+{\rm max} (u,a)\geq
{\rm max} (u_j,a)$ on ${\rm supp}\, \phi\setminus O_\varepsilon$ for all $j$ large enough. Hence by the definition of $C_n$, for all $j$ large enough we have
$$C_n\bigl\{z\in E;\  u(z)>u_j(z)+3\varepsilon\bigr\}\leq$$$$
 {(n!)^2\delta^n\over
\varepsilon^{n+1}}\,\int\limits_{\Omega}\phi\,\bigl(\varepsilon^{n+2}+{\rm max} (u,a)-{\rm max} (u_j,a) \bigr)   
(dd^cv_j)^n+2\,\varepsilon+\varepsilon\,(n!)^2(\varepsilon^{n+2}-a)(-a-t)^n\sup_\Omega |\phi|$$
$$={(n!)^2\delta^n\over
\varepsilon^{n+1}}\,\int\limits_{\Omega}\phi\,\bigl({\rm max} (u,a)-{\rm max} (u_j,a) \bigr)   
\bigl((dd^c {\rm max} (v_j,a+t) )^n-(dd^c{\rm max} ( v,a+t))^n\bigr)$$
$$+{(n!)^2\delta^n\over
\varepsilon^{n+1}}\,\int\limits_{\Omega}\phi\,\bigl({\rm max} (u,a)-{\rm max} (u_j,a) \bigr)   
(dd^c{\rm max} ( v,a+t))^n+{\rm O}(\varepsilon)\longrightarrow {\rm O}(\varepsilon),\ {\rm as\ }j\to\infty,$$
where the last limit follows from Theorem 1 and Corollary 1 in [X1] or [C2]. By the arbitrarility of $\varepsilon> 0$ we get that
$u_j\to u $ in $C_n$ on $E$ as $j\to\infty$, which concludes the proof of Theorem 2.
\bigskip
\noindent {\bf Corollary 1.}  \it  A positive measure $\mu$ in $\Omega$ can be written as $\mu=(dd^cv)^n$ for  
 $v\in {\cal F}(\Omega)$ if and only if $$\mu=(dd^cu_1)^n+ \chi_{\{u_2=-\infty\}}\,(dd^cu_2)^n$$ 
for  $u_1\in {\cal F}^a(\Omega)$ and $ u_2\in {\cal F}(\Omega)$. \rm 
\bigskip
\noindent{\it Proof.} To prove the ``only if" part, by [C2][K] there exists a decreasing
sequence $g_k\in {\cal E}_0(\Omega)$ such that $g_k\geq v$ in $\Omega$ and 
$(dd^cg_k)^n=\chi_{\{v>-k\}}\,(dd^cv)^n$.
Then $u_1:=\lim\limits_{k\to\infty}g_k\in {\cal F}^a(\Omega)$ and 
 $(dd^cu_1)^n=\chi_{\{v\not=-\infty\}}\,(dd^cv)^n$. Hence we have $\mu=(dd^cu_1)^n+ \chi_{\{v=-\infty\}}\,(dd^cv)^n$ .
To prove the ``if" part, From Theorem 1 it turns out that there exists $h\in {\cal F}(\Omega)$ such that $\mu=(dd^cu_1)^n+ (dd^ch)^n.$ By Theorem 5.11 in [C1] there exist a function $\psi\in {\cal E}_0(\Omega)$ and a function $f\in L_{loc}((dd^c\psi)^n)$ such that $(dd^cu_1)^n=f\,(dd^c\psi)^n$. Take a sequence $h_j\in {\cal E}_0(\Omega)$ such that $h_j\searrow h$, $j\to\infty$. Since $\min (f,k^n)\,(dd^c\psi)^n+(dd^ch_j)^n\leq \bigl(dd^c(k\psi+h_j)\bigr)^n $, by [C2][K] there exist $v_j^k\in {\cal E}_0(\Omega)$ such that $(dd^cv_j^k)^n=\min (f,k^n)\,(dd^c\psi)^n+(dd^ch_j)^n$ and hence the comparison theorems in [BT][C1] imply that $0>v^k_j \geq k\,\psi+h\geq u_1+h$. Repeating the proof of Theorem 2 we obtain an increasing sequence $v^k$ in ${\cal F}(\Omega)$ such that  $(dd^cv^k)^n=\min (f,k^n)\,(dd^c\psi)^n+(dd^ch)^n$ and 
$0>v^k\geq u_1+h$. Therefore,  $v:=\bigl(\lim\limits_{k\to\infty}v^k\bigr)^*\in {\cal F}(\Omega)$ and  $\mu=(dd^cv)^n$. The proof of Corollary 1 is complete.
\bigskip
\noindent {\bf Corollary 2.}  \it  For any set $B=\{z_1,z_2,\dots,z_m\}$ of points in $\Omega$ and nonnegative constants $c_1,c_2,\dots,c_m$ there exists a function $u\in PSH(\Omega)\cap L_{loc}^\infty(\Omega\setminus B)$ such that $u=0$ on $\partial\Omega$ and $(dd^cu)^n=\sum\limits_{j=1}^mc_j\,\delta_{z_j}$ in $\Omega$, where $\delta_{z_j}$ denotes the Dirac measure at the point $z_j$.
\rm 
\bigskip
\noindent{\it Proof.} Take the pluricomplex Green function $g_{z_j}$ of $\Omega$ with logarithmic pole at $z_j$ and set $v=\sum\limits_{j=1}^mc_j^{1/n}\,g_{z_j}$. Then $v\in {\cal F}(\Omega)\cap L_{loc}^\infty(\Omega\setminus B)$ and $v=0$ on $\partial\Omega$. By Lemma 5 in [X3] we have that $(dd^cv)^n$ has zero mass at  any point $z\not\in B$ and has mass $c_j\,\delta_{z_j}$ at $z_j$. Therefore, by Theorem 1 we get the required function $u$ and the proof is complete.

\bigskip \bigskip
\bigskip \centerline{\bf References } \bigskip
\bigskip

\noindent [BT] $\,$E.Bedford and B.A.Taylor, {\it A new capacity for plurisubharmonic
functions}. Acta 

Math., {\bf 149} (1982), 1-40.

\noindent [C1] $\,$U.Cegrell, {\it The general definition of the complex Monge-Amp\`ere operator}.
Ann. Inst. 

Fourier {\bf 54}, 159-197 (2004). 

\noindent [C2] $\,$U.Cegrell, {\it
Pluricomplex energy}. Acta Math. {\bf 180:2} (1998), 187-217.

\noindent [K] $\,$ S.Kolodziej, {\it The range of the complex Monge-Amp\`ere operator, II.} Indiana
Univ. 

Math. J. {\bf 44} (1995), 765-782.

\noindent [X1] $\,$Y.Xing, {\it Convergence in Capacity.} Ume\aa\ university, Research Reports No 1,
2007.

\noindent [X2] $\,$Y.Xing, {\it Continuity of the complex Monge-Amp\`ere operator.} 
 Proc. of Amer. Math.  

Soc., {\bf 124} (1996), 457-467.

\noindent [X3] $\,$Y.Xing, {\it The complex Monge-Amp\`ere equations with a countable number of
singular 

points.} Indiana Univ. Math. J, {\bf 48} (1999), 749 - 765.

\bigskip \smallskip
\noindent Department of Mathematics, University of Ume\aa, S-901 87 Ume\aa, Sweden
\smallskip
\noindent E-mail address:\enskip Yang.Xing@mathdept.umu.se

\vfill\eject

\end